%% file: SylowU0-long.tex
\documentclass[a4paper]{article}

\usepackage{amsmath,amssymb,amsfonts,amsthm} 
\usepackage[breaklinks=true]{hyperref}
\usepackage{xspace,enumerate}
\usepackage{burdges}
\usepackage{calc}

\newcommand\Uir{U_{0,r}}
\renewcommand\mathperiod{}

\begin{document}

\title{Sylow theory for $p=0$ in\\solvable groups of finite Morley rank}
\author{Jeffrey Burdges
\thanks{Supported by National Science Foundation grant DMS-0100794
 and Deutsche Forschungsgemeinschaft grant Te 242/3-1.}\\
Fakult\"at f\"ur Mathematik \\ Universit\"at Bielefeld \\ D-33501 Bielefeld, Germany \\
{\tt burdges@math.rutgers.edu}
}

\maketitle

The algebraicity conjecture for simple groups of finite Morley rank,
also known as the Cherlin-Zilber conjecture, states that simple groups
of finite Morley rank are simple algebraic groups over algebraically
closed fields.  In the last 15 years, the main line of attack on this problem
has been Borovik's program of transferring methods from finite group theory.
Borovik's program has led to considerable progress; however, the
conjecture itself remains decidedly open.
In Borovik's program, groups of finite Morley rank are divided into
four types, odd, even, mixed, and degenerate, according to the structure
of their Sylow 2-subgroup.  For {\em even} and {\em mixed type}
the algebraicity conjecture has been proven.

The present paper provides a collection of tools which play a role in
the analysis of {\em odd type} groups, and may have applications in
{\em degenerate type}.  These tools involve the ``0-unipotence''
technology introduced and applied in \cite{Bu03,BuPhd}.
In \cite{Bu03}, the 0-unipotence theory is restricted to those results
needed for the applications in signalizer functor theory.  The tools
presented below develop the general theory further.  These results will
be applied in \cite{Bu05b}, which deals with the Bender method in
minimal simple groups, and again in \cite{BCJ}, where minimal
simple groups of odd type are treated.

A central theme of the present paper is ``Sylow theory for $p=0$'',
in solvable groups of finite Morley rank.  As we will see below, it
would be more accurate to say $p=(0,r)$ here, where the parameter
$r$ represents the ``reduced rank'' of \cite{Bu03} (cf. \S\ref{sec:Uz}).
So, from our point of view, the case $p=0$ splits into infinitely many ``primes''.

The paper is organized as follows.
Section \ref{sec:Uz} recalls the definition of the ``$p$-unipotent radical''
$U_p(G)$ and its analogs for $p=0$, the operators $U_{0,r}$, as introduced
in \cite{Bu03} and \cite{BuPhd}.  Section \ref{sec:Basic} shows
that, in nilpotent groups, several basic properties of the connected
component operator also hold for the $U_{0,r}$ operator, including a
variant of the normalizer condition (Lemma \ref{Unormalizer} below).

Section \ref{sec:NilDecomp} provides the first substantial result,
 a decomposition theorem for nilpotent groups.

\begin{namedtheorem}{Theorem \ref{nildecomp0}}
Let $G$ be a divisible nilpotent group of finite Morley rank,
and let $T$ be the torsion subgroup of $G$.  Then
$$G = d(T) * U_{0,1}(G) * U_{0,2}(G) * \cdots * U_{0,\rk(G)}(G) $$
\end{namedtheorem}

\noindent  One may view this result as an analog of the fact that
a finite nilpotent group is the direct product of its Sylow subgroups.
In section \ref{sec:Criterion}, this decomposition is applied to
produce a new proof of a theorem from \cite{FrJa04}
(see Lemma \ref{nilpotencepre2} below).

In contrast with $p$-groups, connected subgroups of a nilpotent
$\Uir$-group ($G = U_{0,r}(G)$) need not be $\Uir$-subgroup themselves.
Nonetheless, we show that, if a nilpotent $\Uir$-group is generated by
a family of definable subgroups, then it is generated by their $\Uir$-parts 
(Theorem \ref{Ugeneration}).

In the final \S\ref{sec:UzSylow}, we develop the theory of
{\em Sylow $\Uir$-subgroup}.  Here a Sylow $\Uir$-subgroup is
 defined as a maximal nilpotent $\Uir$-subgroup.

\begin{namedtheorem}{Theorem \ref{UirCarter_conj}}
Let $H$ be a connected solvable group of finite Morley rank.
Then the Sylow $\Uir$-subgroups of $H$ are $H$-conjugate.
\end{namedtheorem}

We recall that a nilpotent and self-normalizing subgroup of a group
of finite Morley rank is called a {\em Carter subgroup.}
The proof of Theorem \ref{UirCarter_conj} is parallel to
Wagner's proof of conjugacy for Carter subgroups \cite{Wag}.
Indeed, a Sylow $\Uir$-subgroup is a nilpotent $\Uir$-subgroup
which is self-normalizing as a $\Uir$-group,  in the sense of
Lemma \ref{USylow}.  The Sylow $\Uir$-subgroup theory
is connected to Carter subgroup theory as follows.

\begin{namedtheorem}{Theorem \ref{USylow_decomp}}
Let $H$ be a connected solvable group of finite Morley rank and
let $Q$ be a Carter subgroup of $H$.
Then $\Uir(H') \Uir(Q)$ is a Sylow $\Uir$-subgroup of $H$, and
every Sylow $\Uir$-subgroup has this form for some Carter subgroup $Q$.
\end{namedtheorem}

In \cite{BCJ}, the results above will be used to prove the following.

\begin{namedtheorem}{Theorem}
Let $G$ be a minimal connected simple group of finite Morley rank
 and odd type with a strongly embedded subgroup.
Then $G$ has \Prufer rank one.
\end{namedtheorem}

The results of the present paper will also be used extensively in \cite{Bu05b},
 which will itself be applied in \cite{BCJ}.

\section{Unipotence}\label{sec:Uz}

While there is no intrinsic definition of unipotence in
 a group of finite Morley rank,
there are various analogs of the ``unipotent radical'':
the Fitting subgroup, the $p$-unipotent operators $U_p$, for $p$ prime,
 and their analogs $\Uir$ from \cite{Bu03,BuPhd}.
We recall their definitions.

The {\em Fitting subgroup} $F(G)$ of a group $G$
 of finite Morley rank is the subgroup generated
 by all its nilpotent normal subgroups.
The Fitting subgroup is itself nilpotent and definable
 \cite[Theorem 7.3]{Bel87,Ne91,BN}.
In some contexts, the Fitting subgroup is a reasonable notion
 of unipotence.

A definable subgroup of a connected solvable group $H$
 of finite Morley rank is said to be {\it $p$-unipotent}
if it is a definable connected $p$-group of bounded exponent.

\begin{fact}%
[{\cite[\qCorollary 2.16]{CJ01}; \cite[\qFact 2.36]{ABC97}}]
\label{Upnilpotence} 
Let $H$ be a connected solvable group of finite Morley rank.
Then there is a unique maximal $p$-unipotent subgroup $U_p(H)$ of $H$,
and $U_p(H) \leq F^\o(H)$.
\end{fact}


The present paper is dedicated to the theory of ``characteristic zero''
unipotence introduced in \cite{Bu03}.
We now turn our attention to this (long) definition.

Consider an abelian group $A$ of finite Morley rank.
We say a pair $A_1,A_2 < A$ of proper subgroups of $A$ is
{\it supplemental} if $A_1 + A_2 = A$.
We say that an abelian group of finite Morley rank is
{\it indecomposable} if it has no supplemental pair of
 proper definable subgroups.

\begin{lemma}\label{conind}
An infinite definable indecomposable abelian group is connected.
\end{lemma}

We require the following fact.

\begin{fact}[{\cite[Ex.~11 p.~93 \& Ex.~13c p.~72]{BN}}]\label{nodeftorsion}
Let $G$ be a group of finite Morley rank and let $H \normal G$ be a
definable subgroup.  If $x\in G$ is an element such that $\bar{x}\in G/H$
is a $p$-element, then $x H$ contains a $p$-element.  Furthermore,
if $H$ and $G/H$ are $p^\perp$-groups, then $G$ is a $p^\perp$-group.
\end{fact}

\begin{proof}[Proof of Lemma \ref{conind}]
Let $A$ be a counterexample. 
Since $A$ is disconnected, $A^\o < A$.  Since $A$ is abelian,
there is a definable subgroup $B < A$ which contains $A^\o$ and
has $A/B$ cyclic of prime order $p$.
By Fact \ref{nodeftorsion},
 there is a $p$-element $x\in A$ which lies in a coset of $B$.
So $\gen{x}$ is a definable supplement of $B$ in $A$.
\end{proof}

\begin{fact}[{\cite[\qLemma 2.4]{Bu03}}]\label{decomposition}
Every abelian group of finite Morley rank can be written as a
finite sum of definable indecomposable abelian subgroups.
\end{fact}

The {\it radical} $J(A)$ of a nontrivial definable abelian group of finite
Morley rank is defined to be the maximal proper definable {\em connected}
subgroup without a definable supplement.  The radical $J(A)$ exists and
is unique by \cite[Lemma 2.6]{Bu03}.
In particular, the radical $J(A)$ of a connected indecomposable abelian
group $A$ of finite Morley rank is its unique maximal proper definable
connected subgroup.

We define the {\it reduced (Morley) rank} $\rr(A)$ of a definable abelian
group $A$ to be the Morley rank of the quotient $A/J(A)$,
i.e.\ $\rr(A) = \rk(A/J(A))$.
For a group $G$ of finite Morley rank, and an integer $r$, we define
$$ \Uir(G) = \Genst{A \leq G}{%
\parbox{\widthof{$A$ is a definable indecomposable group,}}%
{$A$ is a definable indecomposable group, \\
\hspace*{10pt} $\rr(A) = r$, and $A/J(A)$ is torsion-free}}\mathperiod $$
We say that $G$ is a {\it $\Uir$-group} (alternatively {\it $(0,r)$-unipotent})
if $U_{0,r}(G)=G$.

The 0-unipotent radical $U_0(G)$ is defined as follows.
Set $\rr_0(G) = \max \{r \mid \Uir(G) \neq 1 \}$, and
 set $U_0(G) = U_{0,\rr_0(G)}(G)$.

We view the reduced rank parameter $r$ as providing a
 {\em scale of unipotence}, with larger values being more unipotent.
By the following fact, the ``most unipotent'' groups are nilpotent.

\begin{fact}%
[{\cite[\qTheorem 2.21]{BuPhd}; \cite[\qTheorem 2.16]{Bu03}}]
\label{nilpotence}
Let $H$ be a connected solvable group of finite Morley rank.
Then $U_0(H) \leq F(H)$.
\end{fact}

%

%

\section{Basic results}\label{sec:Basic}

We begin with some general results on nilpotent $U_{0,r}$-subgroups.
Each result below has a trivial analog for nilpotent $p$-groups, as well as
a less trivial version for {\em connected groups}.

\begin{fact}%
[{\cite[\qLemma 2.12]{BuPhd}; \cite[\qLemma 2.11]{Bu03}}]
\label{Uhom} 
Let $f : G \to H$ be a definable homomorphism between two groups of finite
Morley rank.  Then the following hold.
\begin{conclusions}
\item \textit{(Push-forward)}
  $f(\Uir(G)) \leq \Uir(H)$ is a $\Uir$-group.
\item \textit{(Pull-back)}
  If $\Uir(H) \leq f(G)$ then $f(\Uir(G)) = \Uir(H)$.
\end{conclusions}
In particular, an extension of a $\Uir$-group by a $\Uir$-group is a
$\Uir$-group.
\end{fact}

We may use this fact to quickly show that the lower central series of
a nilpotent $\Uir$-group consists of $\Uir$-groups.

\begin{lemma}[{\cite[\qLemma 2.23]{BuPhd}}]\label{Ucentralseries}
Let $G$ be a nilpotent $\Uir$-group.  Then the subgroups $G^k$
and their quotients $G^k/G^{k+1}$ are $\Uir$-groups for all $k$.
\end{lemma}

It is a standard fact that, in a nilpotent group $G$,
 there is a homomorphism $\ad(g) : G/G' \times G^{k-1}/G^k \to G^k/G^{k+1}$
 induced by the commutator map $(x,y) \mapsto [x,y]$.
Such fact as this, and $[G^i,G^j] \leq G^{i+j+1}$ \cite[Ex.~2d p.~4]{BN},
 are direct consequences of the following identities
 (see also \cite[\qLemma 1.7]{War} and \cite[Ex.~12a p.~6]{BN}).
\begin{eqnarray*} 
[z,uv] = [z,v] [z,u]^v \quad [uv,z] = [u,z]^v [v,z] \\ \relax
[[x,y^{-1}],z]^y [[y,z^{-1}],x]^z [[z,x^{-1}],y]^x = 1
\end{eqnarray*}

\begin{proof}[Proof of Lemma \ref{Ucentralseries}]
We may assume that $G^{k+1}$ is a $\Uir$-group, by downward induction on $k$.
By Lemma \ref{Uhom}, $G/G'$ is a $\Uir$-group.
By Lemma \ref{Uhom}, the image $f_k(G/G',g) \leq G^k/G^{k+1}$
 is a $\Uir$-group for $g\in G^k/G^{k+1}$.
Since these groups generate $G^k/G^{k+1}$,
 the quotient $G^k/G^{k+1}$ is a $\Uir$-group too.
By Lemma \ref{Uhom}, $G^k$ is a $\Uir$-group.
\end{proof}

We also show that the center of a nilpotent group involves all its reduced ranks.

\begin{lemma}\label{Ucenter}
Let $G$ be a nilpotent group of finite Morley rank which satisfies
$\Uir(G) \neq 1$.  Then $\Uir(Z(G)) \neq 1$.
\end{lemma}

\begin{proof}
Suppose that $\Uir(Z_k(G)) = 1$ with $k \geq 1$.
We may assume $k$ is maximal.  Let $H = G/Z_{k-1}(G)$.
By Lemma \ref{Uhom},
 we have $\Uir(Z(H)) = 1$ and $\Uir(Z_2(H)) \neq 1$.
For any $a\in H$,  there is a homomorphism
 $f_a : Z_2(H) \to Z(H)$ given by $x \mapsto [a,x]$.
By Lemma \ref{Uhom},
 $f_a(\Uir(Z_2(H)))$ is a $\Uir$-group,
 so we may assume $f_a(\Uir(Z_2(H))) = 1$.
So $[a,\Uir(Z_2(H))] = 1$ for all $a\in H$,
 and $\Uir(Z_2(H)) \leq Z(H)$, a contradiction.
\end{proof}

We also have a 0-unipotent analog of the connected normalizer condition
\cite[Lemma 6.3]{BN}.

\begin{lemma}[{\cite[\qLemma 2.28]{BuPhd}}]\label{Unormalizer}
Let $G$ be a nilpotent $\Uir$-group.
If $H < G$ is a definable subgroup then $\Uir(N_G(H)/H) > 1$.
\end{lemma}

\begin{proof}
We may assume that $G$ is a counterexample of minimal Morley rank,
and that $H$ is nontrivial.
Then $Z := \Uir(Z(G))$ is nontrivial by Lemma \ref{Ucenter}.
We may assume that $Z \leq H$.  So $K := \Uir(H)$ is nontrivial.
By Lemma \ref{Uhom}, $G/Z$ is a $\Uir$-group.
By minimality, the group
$$ \Uir(N_G(H)/Z) / K  = \Uir(N_{G/Z}(H/Z)) / K $$
is nontrivial.
By Lemma \ref{Uhom}, this is isomorphic to $(\Uir(N_G(H)) / Z) / (K/Z)$,
 which is isomorphic to $\Uir(N_G(H)) / K$, a contradiction.
\end{proof}

Lemma \ref{Uhom} also yields the following generalization of Fact \ref{nilpotence}.

\begin{lemma}\label{nilpotencecor}
Let $H$ be a connected solvable group of finite Morley rank.
Let $r$ be the maximal reduced rank such that $\Uir(H) \not\leq Z_n(H)$
for any $n$.  Then $\Uir(H) \leq F(H)$.
\end{lemma}

\begin{proof}
Suppose that $U_{0,s}(H) \leq Z(H)$ for all $s>r$.
Then $\Uir(H/Z^\o(H))$ is nilpotent by Fact \ref{nilpotence}.
By Lemma \ref{Uhom},
 $$ \Uir(H) / (\Uir(H) \cap Z^\o(H)) \cong
 \Uir(H) Z^\o(H) / Z^\o(H) \cong \Uir(H/Z^\o(H)) \mathperiod $$
So $\Uir(H)$ is a central extension of a nilpotent group,
 and hence nilpotent itself (see Fact \ref{basicnilcrit} below).
\end{proof} 

We recall that the {\em Frattini subgroup} $\Phi(G)$ of
 a group $G$ of finite Morley rank is defined as
 the intersection of all maximal proper {\em definable connected}
subgroups of $G$ (see also \cite[Definition 5.6]{Fre00a}).
This definition is designed to yield the following fact.

\begin{fact}[{\cite[\qLemma 5.7]{Fre00a}}]\label{frattini_gen}
Let $G$ be a connected group of finite Morley rank.
Suppose $H$ is a definable subgroup of $G$ such that $G = H \Phi(G)$.
Then $H = G$.
\end{fact}


\begin{corollary}\label{Phi_vs_J}
Let $A$ be a connected abelian group of finite Morley rank.
Then $\Phi^\o(A) = J(A)$.
\end{corollary}

\begin{proof}
$\Phi^\o(A)$ has no definable supplement, by Fact \ref{frattini_gen},
 so $\Phi^\o(A) \leq J(A)$.
Conversely, consider a maximal proper definable connected
subgroup $M$ of $A$.  Then $A/M$ is a minimal group.
Since $J(A)$ has no supplement in $A$, $J(A) \leq M$.
Since $J(A)$ is connected, $J(A) \leq \Phi^\o(A)$.
\end{proof}

Olivier \Frecon has pointed out that nilpotent $\Uir$-groups
 are ``homogeneous'' modulo their Frattini subgroup.

\begin{lemma}[cf.\ \Frecon]\label{Phi_homogeneity}
Let $G$ be a nilpotent $\Uir$-group of finite Morley rank.
Then $U_{0,s}(G) \leq \Phi(G)$ for $s \neq r$,
 and $\Phi(G)$ contains all torsion elements of $G$.
\end{lemma}

\begin{proof}
Otherwise,
there is a maximal proper definable connected subgroup $M$ of $G$ such that
 either $G/M$ contains torsion, by Fact \ref{nodeftorsion},
 or $U_{0,s}(G/M) \neq 1$ for some $s \neq r$, by Fact \ref{Uhom}.
But $G/M$ is a minimal group, so $\Uir(G/M) = 1$ in either case,
 a contradiction to Fact \ref{Uhom}.
\end{proof}

We note that Facts \ref{frattini_gen} and \ref{Uhom}
 provide a converse to Lemma \ref{Phi_homogeneity}.
i.e.\ $G$ is a $\Uir$-group if $G/\Phi(G)$ is a $\Uir$-group.

\Frecon's lemma provides an easy proof of the following important
generation theorem.  Our original proof of this result required a
Krull-Schmidt theorem. 

\begin{theorem}[{\cite[\qTheorem 2.41]{BuPhd}}]\label{Ugeneration}
Let $G$ be a nilpotent $\Uir$-group.  Let $\{ H_i \}_{i \leq n}$ be a
family of definable subgroups such that $G = \gen{ \cup_i H_i }$.
Then $G = \gen{ \cup_i \Uir(H_i) }$.
\end{theorem}

\begin{proof} 
Consider the quotient $\bar{G} = G / \Phi(G)$.  Then $\bar{G}$
is generated by the groups $\bar{H_i} = H_i \Phi(G)/\Phi(G)$.
Since $G$ is a $\Uir$-group, Lemma \ref{Phi_homogeneity} says that
 $\bar{G}$ is torsion-free and $U_{0,s}(\bar{G}) = 1$ for $s \neq r$,
 by Fact \ref{nodeftorsion} and Lemma \ref{Uhom}, respectively.
So $\Uir(\bar{H_i}) = \bar{H_i}$.
By Fact \ref{Uhom}, $G = \gen{ \cup_i \Uir(H_i) } \Phi(G)$.
So the result follows from Fact \ref{frattini_gen}.
\end{proof}

\section{A nilpotent structure theorem}\label{sec:NilDecomp}

We now use the 0-unipotence theory to find a canonical decomposition
of a connected nilpotent group of finite Morley rank.
This generalizes the fact that a finite nilpotent group is the product
of its Sylow $p$-subgroups, a result which follows from the next fact.

\begin{fact}\label{Upqcommutator}
Let $N$ be a nilpotent group, let $p,q$ be distinct primes,
 let $P$ be a $p$-subgroup of $N$, and let $Q$ be a $q$-subgroup of $N$.
Then $[P,Q]=1$.
\end{fact}

As suggested by the analogy with prime numbers,
 a similar result holds for $U_{0,r}$-groups of different reduced ranks $r$.

In the proof, we will use the following fact.

\begin{fact}[{\cite[\qCorollary 5.29]{BN}; \cite{Zil77}}]\label{concommutator}
Let $H$ be a definable connected subgroup of a group $G$ of finite Morley rank
and let $X \subset G$ be any subset of $G$.  Then $[H,X]$ is definable and
connected.
\end{fact}

\begin{lemma}[{\cite[\qLemma 2.30]{BuPhd}}]\label{Uzrcommutator}
Let $N$ be a nilpotent group of finite Morley rank.
Let $A,B \leq N$ be definable connected indecomposable abelian groups
with $\rr(A) \neq \rr(B)$.  Then $[A,B] = 1$.
\end{lemma}

To prove this, we recall a technical version of Fact \ref{Uhom}-1.

\begin{fact}%
[{\cite[\qLemma 2.9]{BuPhd}; \cite[\qLemma 2.9]{Bu03}}]
\label{Indpushforward} 
Let $f : A \to G$ be a definable homomorphism between two definable connected
groups inside a structure of finite Morley rank.  Suppose that $A$ is abelian
and indecomposable.  Then $f(A)$ is indecomposable and $f(J(A)) = J(f(A))$.
\end{fact}

\begin{proof}[Proof of Lemma \ref{Uzrcommutator}]
Suppose toward a contradiction that $[A,B] \neq 1$.
Let $k$ be the maximal integer such that $[A,B]\leq N^k$.
We may assume that $N^{k+1} = 1$ by Fact \ref{Indpushforward},
 and then $[A,B]$ is central in $N$.
By Fact \ref{concommutator}, the group $[A,B]$ is connected and definable.
Let $M < [A,B]$ be a proper connected subgroup of maximal Morley rank.
We may also assume $M = 1$, by Fact \ref{Indpushforward},
 since $M \normal N$. 
So $[A,B]$ is indecomposable and $J([A,B]) = 1$.
There are $a\in A$ and $b\in B$ such that $[a,b] \neq 1$.
There are homomorphisms $f : A \to [A,B]$ and $g : B \to [A,B]$
 given by $x \mapsto [x,b]$ and $y \mapsto [a,y]$, respectively.
Since $[a,b] \neq 1$ and $A,B$ are connected,
 $f(A)$ and $g(B)$ are nontrivial connected subgroups of $[A,B]$.
Since $[A,B]$ is minimal, $f(A) = [A,B] = g(B)$.
Since $[A,B]$ is indecomposable, Lemma \ref{Uhom} says
$\rr(A) = \rr([A,B]) = \rr(B)$, a contradiction.
\end{proof}

\begin{theorem}\label{nildecomp0}
Let $G$ be a divisible nilpotent group of finite Morley rank,
and let $T$ be the torsion subgroup of $G$.  Then
$$G = d(T) * U_{0,1}(G) * U_{0,2}(G) * \cdots * U_{0,\rk(G)}(G) $$
\end{theorem}

\begin{proof}
By \cite[Theorem 6.9]{Ne91nilp,BN}, $T$ is central in $G$,
 and thus $d(T) \leq Z(G)$.
Since $U_{0,r}(G)$ and $U_{0,s}(G)$ commute whenever $r \neq s$
by Lemma \ref{Uzrcommutator}, we have a central product
$$ D := d(T) * U_{0,1}(G) * U_{0,2}(G) * \cdots * U_{0,\rk(G)}(G) $$
and we claim $D = G$.
Since $d(T) \normaleq G$ and $\Uir(G) \normal G$,
 we have $D \normal G$.
Suppose toward a contradiction that $D < G$.
Since $G/D$ is torsion-free by Fact \ref{nodeftorsion},
$U_0(G/D) \neq 1$ and $U_{0,\rr_0(G/D)}(G) \not\leq D$
by Lemma \ref{Uhom}, a contradiction.
\end{proof}

We combine Theorem \ref{nildecomp0} and Fact \ref{Upqcommutator}
 as follows:

\begin{corollary}[{\cite[\qTheorem 2.31]{BuPhd}}]\label{nildecomp}
Let $G$ be a nilpotent group of finite Morley rank.
Then $G = D * B$ is a central product of definable characteristic
 subgroups $D,B \leq G$ where $D$ is divisible and
 $B$ has bounded exponent (it is connected iff $G$ is connected).
Let $T$ be the torsion part of $D$.
Then we have decompositions of $D$ and $B$ as follows.
\begin{eqnarray*}
D &=& d(T) * U_{0,1}(G) * U_{0,2}(G) * \cdots \\
B &=& U_2(G) \times U_3(G) \times U_5(G) \times \cdots
\end{eqnarray*}
\end{corollary}

\begin{proof}
The group $G$ is a central product of groups $D$ and $B$
 by \cite[Theorem 6.8, Ex.\ 7 p.~71, 12d p.~6]{BN}.
Theorem \ref{nildecomp0} provides the decomposition of $D$.
By Fact \ref{Upqcommutator},
 $C := U_2(G) \times U_3(G) \times U_5(G) \times \cdots$ is a direct sum.
Since $U_p(G) \normal G$ for prime $p$, $C \normal G$.
By Fact \ref{nodeftorsion}, $U_p(G/C) = 1$ for any $p$ prime.
So $B = C$.
\end{proof}

As a corollary, we have a version of Fact \ref{concommutator}
for normal nilpotent $\Uir$-groups.  As in \S\ref{sec:Basic},
the operator $U_{0,r}$ replaces the connected component operator.
Such results are trivial for $U_p$ because a connected
subgroup of a $p$-unipotent group is $p$-unipotent.

\begin{corollary}\label{Unilcommutator}
Let $G$ be a solvable group of finite Morley rank, let $S \subset G$ be
any subset, and let $H$ be a nilpotent $\Uir$-group which is normal in $G$.
Then $[H,S] \leq H$ is a $\Uir$-group.
\end{corollary}

\begin{proof}
We may assume that $S = \{s\}$ consists of a single element.
Let $Y := \Uir([H,s])$ be the largest $\Uir$-subgroup of $[H,s]$.
We may assume that $H' \leq Y$ by induction on the nilpotency class of $H$,
 and so $Y \normal H$. 
Consider the commutator map $f : H \to [H,s]/Y$ given by $x \mapsto [x,s]$.
Since $H \normal G$, the group $[H,s] \leq H$ is nilpotent.
By Theorem \ref{nildecomp},
 $[H,s] = d(T) * U_{0,1}([H,s]) * U_{0,2}([H,s) * \cdots$,
 where $T$ is the torsion part of $[H,s]$.
Since $H$ is a $\Uir$-group,
 $U_{0,k}([H,s]) \leq U_{0,k}(H) \leq Z(H)$ for $k \neq r$,
 by Theorem \ref{nildecomp} again.  Similarly, $d(T) \leq Z(H)$.
So, for any $u\in H$, there is a $y\in Y$ and a $z\in [H,s] \cap Z(H)$
 such that $[u,s] = z y$,
i.e.\ $z$ is the product of elements from the other factors of $[H,s]$.
For any $v\in H$,
$$ [uv,s] = [u,s]^v [v,s] = (zy)^v [v,s] = z y^v [v,s] \in [u,s] [v,s] Y\mathperiod $$
So the map $f$ is a homomorphism.
Since $H$ is a $\Uir$-group,
 the image $[H,s]/Y$ is a $\Uir$-group by Fact \ref{Uhom}.
Since $Y$ is also a $\Uir$-group,
 $[H,s]$ is a $\Uir$-group by Fact \ref{Uhom}.
\end{proof}

\section{A nilpotence criterion}\label{sec:Criterion}

In \cite{FrJa04}, \Frecon and Jaligot used the 0-unipotent
radical to build Carter subgroups of non-solvable groups of
finite Morley rank.  Their idea was to build an almost self-normalizing
nilpotent group in stages by monitoring the reduced rank parameter,
first taking the least unipotent groups, and then adding increasingly
unipotent groups which normalize them.
This scale can be summarized as follows:
\begin{enumerate}
\item a divisible torsion group $S_i$ is a ``non-unipotent'' group
  because it is never acted upon by a connected group.
\item the reduced rank parameter $r$ measures the degree of
  unipotence for nilpotent $\Uir$-groups.
\item a $p$-unipotent group is always unipotent.
\end{enumerate}
These points are represented in \cite[Corollary 3.8]{FrJa04},
but the proof here clarifies matters by employing Theorem \ref{nildecomp0}.
Lemma \ref{nilpotencepre2}, which justifies the second point's
scale of 0-unipotence, is the only result of this section which will
play a role in the remainder of the paper.

\begin{proposition}\label{nilpotence2}
Let $G = \gen{ G_1,\ldots,G_n }$ be a group of finite Morley rank which is
generated by connected definable nilpotent groups $G_i$, each of which satisfies
one of the following.
\begin{hypotheses}
\item[(i)] $G_i = d(S_i)$ for some divisible torsion group $S_i$,
\item[(ii)]  $U_{0,r_i}(G_i) = G_i$ for some natural number $r_i$, or
\item[(iii)]  $U_{p_i}(G_i) = G_i$ for some $p_i$ prime.
\end{hypotheses}
Suppose the $G_i$ are ordered from least unipotent to most unipotent,
i.e.\ groups of type (i) come before groups of type (ii) or (iii),
 and the groups of type (ii) have nondecreasing reduced ranks,
i.e.\ $r_i \geq r_j$ whenever $i > j$.
Suppose further that $G_i$ normalizes $G_j$ whenever $i > j$.

Then $G$ is nilpotent.
\end{proposition}

Fact \ref{Upnilpotence} suggests that $p$-unipotent groups, type (iii),
are more unipotent then either type (i) or type (ii).  However,
a field of finite Morley rank of positive characteristic has no definable
torsion-free section of its multiplicative group \cite{Wag01}.
From this, it will follow that a nilpotent $\Uir$-group never acts
on a $p$-unipotent group.
As a result, one may allow groups of type (iii) to precede groups
 of type (ii).


The proposition will use standard tricks for building a central series.

\begin{fact}[{\cite[\qFact 2.5]{FrJa04}}]\label{basicnilcrit}
Let $G = HL$ be a group with $H$ and $L$ nilpotent subgroups and $H \normal G$.
Let $1 = H_0 \leq H_1 \leq \cdots \leq H_n = H$ be a series of $G$-normal
subgroups of $H$.  Suppose that the groups $(H_i/H_{i-1}) \rtimes L$ are
nilpotent for $1 \leq i \leq n$.  Then $G$ is nilpotent.
\end{fact}  

We recall that connected groups do not act upon divisible torsion,
 to dispense with type (i).

\begin{fact}[{\cite[Ex. 3 p.\ 148]{BN}}]\label{cent_divtorsion}
Let $G$ be a connected group of finite Morley rank. 
Then the divisible torsion part of $F(G)$ is central in $G$.
\end{fact}

We capture the behavior of type (ii) as follows.

\begin{lemma}[{compare \cite[\qProposition 3.7]{FrJa04}}]\label{nilpotencepre2}
Let $G = H T$ be a group of finite Morley rank with
 $H \normal G$ a nilpotent $U_{0,r}$-group and
  $T$ a nilpotent $U_{0,s}$-group for some $s \geq r$.
Then $G$ is nilpotent.
\end{lemma}

\begin{proof}
Since $[H,T]$ is a $\Uir$-group by Corollary \ref{Unilcommutator},
 we may suppose inductively that $[H,T] \rtimes T$ is nilpotent.
So we may assume that $[H,T] = H$,
 as the lemma follows by factoring out $[H,T]$ otherwise.
Let $\bar{H} = H/\Phi(H)$.
So $s \geq \rr_0(\bar{H})$ by Lemma \ref{Phi_homogeneity}.
As $T \leq F(\bar{H} \rtimes T)$ by Fact \ref{nilpotence}, 
 $\bar{H} \rtimes T$ is nilpotent.
So $[\bar{H},T] = \bar{H}$ contradicts Theorem \ref{nildecomp0}.
\end{proof}

%
%

We finally treat the relationship between types (ii) and (iii).

\begin{lemma}\label{Uz_on_Up}
Let $G$ be a connected solvable group of finite Morley rank.
Suppose that $S$ is a nilpotent $\Uir$-subgroup of $G$, and
 that $G = U_p(G) S$ for some $p$ prime.
Then $G$ is nilpotent, and $[U_p(G),S]=1$.
\end{lemma}

\begin{proof}
There is a normal series $\{A_i\}_{i=0}^n$ for $U_p(G)$
 such that $A_i/A_{i-1}$ is {\em $S$-minimal},
 i.e.\ minimal subject to being definable, infinite, and $S$-normal.
Suppose toward a contradiction that $C_S(A_i/A_{i-1}) < S$.
By the Zilber Field Theorem \cite[Theorem 9.1]{BN},
 there is an algebraically closed field $k$ with $k_+ \cong A_i/A_{i-1}$,
 and $S/C_S(A_i/A_{i-1}) \hookrightarrow k^*$.
So $k$ is a definable field of characteristic $p$,
 with a definable torsion-free section of $k^*$,
 in contradiction to \cite[Corollary 9]{Wag01}.
So $G = U_p(G) S$ is nilpotent by Fact \ref{basicnilcrit}.
By Corollary \ref{nildecomp}, $[U_p(G),S]=1$.
\end{proof}

\smallskip

\begin{proof}[Proof of Proposition \ref{nilpotence2}]
Let $N = \gen{ G_1,\ldots,G_{n-1} }$.
We assume inductively that $N$ is nilpotent.
So $G_n$ normalizes $N$ and $G = \gen{ G_n \cup N }$ is solvable.

Suppose first that $G_n$ is of type (i), i.e.\ $G_n = d(S_n)$ 
for some divisible torsion group $S_{n-1}$.
Then $G_i$ for $i<n$ are all of type (i) by our assumption;
and thus $G_i \leq Z(G)$ for $i<n$ by Fact \ref{cent_divtorsion}.
So $N\leq Z(G)$ and $G$ is abelian.

For the type (ii) case, suppose $U_{0,r_n}(G_n) = G_n$ for some $r_n\in\N$.
Let $H_i := \gen{G_j : j \leq i} / \gen{G_j : j<i}$ for $i<n$.
As above, $H_i$ is abelian whenever $G_i$ is of type (i).
By Fact \ref{Uhom}, the group $H_i$ is a $U_{0,r_i}$-group whenever
 $G_i$ is of type (ii), i.e.\ $G_i$ is a $U_{0,r_i}$-group.
Clearly $H_i$ is $p_i$-unipotent when $G_i$ is of type (iii) too.
We observe that $G_n$ acts on $H_i$ for $i<n$.
By Fact \ref{cent_divtorsion}, Lemma \ref{nilpotencepre2},
 or Lemma \ref{Uz_on_Up} (as appropriate)
 $H_i \rtimes G_n$ is nilpotent for $i<n$.
So $G$ is nilpotent by Fact \ref{basicnilcrit}.

For the type (iii) case, suppose $U_{p_n}(G_n) = G_n$ for some $p_n$ prime.
The result follows immediately from the fact that $G_n\leq F(G)$
by Fact \ref{Upnilpotence}.
\end{proof}

A convenient corollary is provided by Corollary \ref{nildecomp}.

\begin{corollary}[{compare \cite[\qProposition 3.7]{FrJa04}}]
Let $G = H T$ be a group of finite Morley rank,
 with $H$ and $T$ definable and nilpotent.
Suppose that $T$ is generated by the groups
 $\Uir(T)$ with $r \geq \rr_0(H)$, and $U_p(T)$ with $p$ prime.
Then $G$ is nilpotent.
\end{corollary}

\begin{proof}
Corollary \ref{nildecomp} provides $H$ and $T$ with decompositions
 $H_1,\ldots,H_n$ and $T_1,\ldots,T_m$ which are compatible with
 the hypotheses of Proposition \ref{nilpotence2}.
Clearly $T$ normalizes each factor of $H$.
So $G = \gen{H_1,\ldots,H_n,T_1,\ldots,T_m}$ is nilpotent
 by Proposition \ref{nilpotence2}.
\end{proof}

\section{Sylow $\Uir$-subgroups}\label{sec:UzSylow}

We now turn to our theory of Sylow $\Uir$-subgroups.

\begin{definition}
Let $G$ be a group of finite Morley rank.  We call a maximal definable 
nilpotent $\Uir$-subgroup a {\em Sylow $\Uir$-subgroup} of $G$.
\end{definition}

Sylow $\Uir$-subgroup are an analog of Carter subgroups in the following sense.

\begin{lemma}[{\cite[\qLemma 4.18]{BuPhd}}]\label{USylow}
Let $H$ be a group of finite Morley rank.
Then the Sylow $\Uir$-subgroups of $H$ are exactly those
 nilpotent $\Uir$-subgroups $S$ such that $\Uir(N_G(S)) = S$.
\end{lemma}

\begin{proof}
Let $S$ be a Sylow $\Uir$-subgroup of $H$.
Suppose toward a contradiction that $\Uir(N_H(S)) > S$.
Then there is an abelian $\Uir$-group $A \leq N_H(S)$ with $A \not\leq S$.
The group $S A$ is nilpotent by Lemma \ref{nilpotencepre2},
contradicting the maximality of $S$.
Conversely, a nilpotent $\Uir$-subgroup $S$ with $\Uir(N_G(S)) = S$
 is a Sylow $\Uir$-subgroup by Lemma \ref{Unormalizer}.
\end{proof}

To show that Sylow $\Uir$-subgroups are conjugate in a solvable
groups, we follow Wagner's proof of conjugacy for Carter subgroups \cite{Wag}.

Given a group $G$, we recall that a subgroup $L$ is said to be a
{\em complement} to a normal subgroup $N \normal G$ in $G$
if $N L = G$ and $N \cap L = 1$.

\begin{fact}[{\cite[\qCorollary 5.5.4]{Wag}}]\label{complements}
Let $G$ be a group of finite Morley rank with a definable normal abelian
subgroup $K$ such that $G/K$ is nilpotent.  Suppose some $g\in G$ acts
without fixed points on $K$.  Then there is a nilpotent complement
$C$ in $G$ for $K$, any nilpotent subgroup $L$ with $G = K L$ is a
complement for $K$, and any two complements to $K$ in $G$ are conjugate.
\end{fact}  

\begin{lemma}\label{UirCarter_quotient}
Let $H$ be a connected solvable group of finite Morley rank and
 let $K$ be a definable $H$-minimal subgroup of $Z(F(H))$.
Then the Sylow $\Uir$-subgroups of $H/K$ are exactly the images of
the Sylow $\Uir$-subgroups of $H$.
\end{lemma}

%
%
\begin{proof}
Suppose that $H$ is a counterexample of minimal Morley rank.
We observe that $J(K) = 1$.  Let $\bar{H} = H/K$.
Consider a Sylow $\Uir$-subgroup $S$ of $H$.
We will show, using the criteria of Lemma \ref{USylow}, that
 $\bar{S} = S K/K$ is a Sylow $\Uir$-subgroup of $\bar{H}$.

First, suppose that $K$ is {\em centralized} by $S$.
Then $K \leq N_H(S)$ and $\Uir(S K) = S$.
So $S \normal N_H(S K)$, 
 and $N_H(S \bmod K) = N_H(S K) \leq N_H(S)$.
By Lemma \ref{Uhom},
 $\Uir(N_H(S \bmod K)) \leq \Uir(N_H(S)) = S$ and
$\Uir(N_{\bar{H}}(\bar{S})) = \bar{S}$, as desired.

Next, suppose that $S$ acts {\em nontrivially} on $K$.
Since $K$ is $H$-minimal,
 $K$ is also $T$-minimal where $T := H/C_H(K)$.
By the Zilber Field Theorem \cite[Theorem 9.1]{BN},
there is an algebraically closed field $k$ with $K \cong k_+$,
$T \hookrightarrow k^*$, and our action is the field action.
Therefore any element of $S \setminus C_H(K)$
 acts without fixed points on $K$.
By Fact \ref{complements}, $S$ is a complement to $K$ in $S K$,
 and the complements to $K$ are all conjugate.
So $N_H(S K) = N_H(S) K$ by a Frattini argument.
The quotient map $H \to \bar{H}$ now gives a surjective
 homomorphism $f : N_H(S) \to N_{\bar{H}}(\bar{S})$.
By Lemma \ref{Uhom},
\[ \Uir(N_{\bar{H}}(\bar{S})) = f(\Uir(N_H(S))) = f(S) = \bar{S} \qedhere \]
\end{proof}

\begin{theorem}[{\cite[\qTheorem 4.16]{BuPhd}}]\label{UirCarter_conj}
Let $H$ be a connected solvable group of finite Morley rank.
Then the Sylow $\Uir$-subgroups of $H$ are $H$-conjugate.
\end{theorem}

\begin{proof}
Let $S_1,S_2$ be Sylow $\Uir$-subgroups of $H$ and
 let $K$ be an $H$-minimal subgroup of $Z(F(H))$.
By Lemma \ref{UirCarter_quotient},  the quotients
 $S_1 K/K$ and $S_2 K/K$ are Sylow $\Uir$-subgroups of $H/K$.
So, by induction on $\rk(H)$, we may assume $S_1 K = H = S_2 K$.
If $K$ is centralized by $S_1$, then $K$ is also centralized by $S_2$,
 so $S_1 = \Uir(S_1 K) = \Uir(S_2 K) = S_2$.
Thus we may assume that $S_1$ and $S_2$ act nontrivially on $K$.
As in Lemma \ref{UirCarter_quotient},
the Zilber Field Theorem shows that there are elements
 of both $S_1$ and $S_2$ which act without fixed points on $K$.
By Fact \ref{complements}, $S_1$ and $S_2$ are complements to $K$
 in $S_1 K$, and the complements to $K$ in $S_1 K$ are all conjugate, as desired.
\end{proof}

As usual, we have the Frattini argument.

\begin{corollary}\label{UirCarter_frattini}
Let $H$ be a solvable group of finite Morley rank, let $K \normal H$ be
a connected definable subgroup, and let $S$ be a Sylow $\Uir$-subgroup
of $K$.  Then $H = N_H(S) K$.
\end{corollary}

\begin{proof}
For any $h\in H$, $S^h$ is a Sylow $\Uir$-subgroup of $K$.
By Theorem \ref{UirCarter_conj}, $S^h = S^k$ for some $k\in K$.
So $hk^{-1} \in N_H(S)$ and $H = N_G(S) K$.
\end{proof}

%

We show next that a Sylow $\Uir$-subgroup of a solvable group
is a product of a ``semisimple'' and a ``unipotent'' factor.

\begin{theorem}[{\cite[\qLemma 4.19]{BuPhd}}]\label{USylow_decomp}
Let $H$ be a connected solvable group of finite Morley rank and
let $Q$ be a Carter subgroup of $H$.
Then $\Uir(H') \Uir(Q)$ is a Sylow $\Uir$-subgroup of $H$, and
every Sylow $\Uir$-subgroup has this form for some Carter subgroup $Q$.
\end{theorem}

\begin{proof}
First, the group $D := \Uir(H') \Uir(Q)$ is nilpotent
 by Lemma \ref{nilpotencepre2}.  
Let $K$ be a Sylow $U_{0,r}$-subgroup of $H$ which contains $D$.
Let $\bar{H} := H/H'$ and let $\bar{K} := K H'/H'$.
By Fact \ref{Uhom}, $\bar{K} \cong K / K \cap H'$ is a $\Uir$-group.
Since the quotient map $f : Q \to Q H'/H'$ is surjective,
$\bar{K} \leq \Uir(\bar{H}) = f(\Uir(Q))$ by Fact \ref{Uhom},
and hence $K \leq H' \Uir(Q)$.
So $K = \Uir(Q) (K \cap H')$.
Now $K = \Uir(Q) \Uir(K \cap H') \leq D$ by Theorem \ref{Ugeneration},
 as desired.
The final comment now follows from Theorem \ref{UirCarter_conj}.
\end{proof}

\begin{corollary}\label{USylow_decomp_Greg}
Let $H$ be a connected solvable group of finite Morley rank such that
$\Uir(H') = 1$.  Then the Sylow $\Uir$-subgroups of $H$ are
contained in Carter subgroups of $H$.
\end{corollary}

\begin{corollary}\label{USylow_decomp_Frecon}
Let $H$ be a connected solvable group of finite Morley rank and
let $S$ be a Sylow $\Uir$-subgroup of $H$.
Then $N_H(S)$ contains a Carter subgroup of $H$.
\end{corollary}

\begin{proof}
Let $Q$ be a Carter subgroup of $H$.  Then $Q$ normalizes $\Uir(Q)$
and $\Uir(H')$, so $Q$ normalizes some Sylow $\Uir$-subgroup.  Now
the result follows by Theorem \ref{UirCarter_conj}.
\end{proof}

Olivier \Frecon has pointed out one may prove Theorem \ref{UirCarter_conj},
Theorem \ref{USylow_decomp}, and Corollary \ref{USylow_decomp_Frecon}
in the reverse order from that shown above.
\Frecon observes that a Sylow $\Uir$-subgroup is normal in its
generalized centralizer, by \cite[Corollaire 5.17]{Fre00a}.
So Corollary \ref{USylow_decomp_Frecon} follows from the facts
that the a generalized centralizer of a nilpotent group is abnormal,
\cite[Corollaire 7.4]{Fre00a}, and abnormal subgroups contain
Carter subgroups \cite[Th\'eor\`eme 1.2]{Fre00a}.
The remainder of \Frecon's reversal is left as an exercise to the reader
(see \cite[Corollaire 5.20]{Fre00a}).

In any case, Corollary \ref{USylow_decomp_Frecon}
is analogous to the following fact for ordinary Sylow subgroups.

\begin{fact}[{\cite[\qCorollary 7.15]{Fre00a}}]\label{Carter_Sylow}
Let $H$ be a connected solvable group of finite Morley rank,
 and let $R$ be a Sylow $p$-subgroup of $H$.
Then $N_H(R)$ contains a Carter subgroup of $H$.
\end{fact}

\begin{corollary}\label{USylow_decomp_end}
Let $H$ be a connected solvable group of finite Morley rank such that
$\Uir(H') = 1$.  Then each Sylow $\Uir$-subgroup of $H$
commutes with some Sylow $p$-subgroup of $H$.
\end{corollary}

Olivier \Frecon has pointed out that one may replace Sylow $p$-subgroups
by Hall $\pi$-subgroups in this corollary, by using \cite[4.17, 4.22]{Fre00b}.

\begin{proof}
By \cite[Theorem 9.29 and \S6.4]{BN},
 any Sylow $p$-subgroup $P$ of $H$ has the form $T * U_p(H)$
 where $T$ is a divisible abelian $p$-subgroup of $H$.
By Fact \ref{Carter_Sylow},
 $N_H(P)$ contains a Carter subgroup $Q$ of $H$.
By Fact \ref{cent_divtorsion}, $T \leq Z(Q)$.
By Corollary \ref{USylow_decomp_Greg},
 $\Uir(Q)$ is a Sylow $\Uir$-subgroup of $H$.
Now $\Uir(Q)$ centralizes $P = T U_p(H)$
 by Lemma \ref{Uz_on_Up}, as $U_p(H) \normal H$,  
The corollary now follows from Theorem \ref{UirCarter_conj}.
\end{proof}

The absence of a central product in Theorem \ref{USylow_decomp}
is conspicuous.  It is reasonable to ask if there are subgroups
$A \leq \Uir(F(H))$ and $B \leq \Uir(Q)$ such that $A B$ is a
Sylow $\Uir$-subgroup of $H$, and $[A,B]=1$.
A positive answer to this is linked to the problem of building a true
notion of unipotent radical in a solvable groups of finite Morley rank.
However, an extension of the group $H$ below suggests that the answer
is no, even for solvable groups which are vaguely ``linear.''
$$
  H := \begin{pmatrix} k_* & k_+ & k_+ \\ & 1 & k_+ \\ & & 1 \\ \end{pmatrix}
\quad
  Q := \begin{pmatrix} k_* & 0 & 0 \\ & 1 & k_+ \\ & & 1 \\ \end{pmatrix}
$$
In this group, $U_{0,1}(H') \cong (k_+)^2$ and $U_{0,1}(Q) \cong k_* \times k_+$
are both abelian, but $H$ still has a nonabelian Sylow $\Uir$-subgroup.
Of course, one could still take $A = F(H)$ and $B=1$.
Now consider a field $l$ with $k_+ \hookrightarrow l_*$.
We can envision extending $H$ by a copy of $l_+$, and asking that the
copy of $k_+$ inside $Q$ act as a group of field automorphisms on $l_+$.
The Fitting subgroup of this extension of $H$ would be abelian, as would
the Carter subgroup $Q$, but the Sylow $U_{0,1}$-subgroup would be
nonabelian.

In practice, Sylow $\Uir$-subgroups are useful because one contains
 any $\Uir$-subgroup, but they still behave like Carter subgroups,
and Theorem \ref{USylow_decomp} links to the usual Carter subgroup.
In \cite{BuPhd}, the method of \cite{FrJa04} is used to construct
more general Carter $U_{\pi,\rho,\tau}$-groups, with $\pi$ and $\tau$
sets of primes and $\rho$ a set of reduced ranks, corresponding to the
three clauses of Proposition \ref{nilpotence2}. 
For $\tau = \pi = \emptyset$, and $\rho = \{r\}$,
these reduce to our Sylow $\Uir$-subgroups.
In \cite{BuPhd}, a conjugacy theorem is proven for such groups.
As Theorem \ref{USylow_decomp} and its corollaries do not generalize,
Carter $U_{\pi,\rho,\tau}$-groups have not yet found applications.

\input{SylowU0-6-KS}
% ...

\section*{Acknowledgments}

I thank my advisor Gregory Cherlin for direction during my thesis work,
of which this article is an outgrowth, and for guidance during its later
developments, as well as for his collaboration in \cite{BCJ}, where
the theory of Sylow $\Uir$-subgroups is applied.

I thank Tuna \Altinel and Olivier \Frecon for their attentive readings,
corrections, and valuable suggestions.  I also thank Eric Jaligot and
Alexandre Borovik for their encouragement and comments.

Financial support for this work comes from NSF grant DMS-0100794,
DFG grant Te 242/3-1, and the Isaac Newton Institute, Cambridge.
I also thank the following institutions for their hospitality:
University of Birmingham, University of Manchester,
IGD at Universit\'e Lyon I, Isaac Newton Institute, Cambridge,
and CIRM at Luminy.

\small
\bibliographystyle{alpha} 
\bibliography{burdges,fMr}


\end{document}

%% file: SylowU0-6-KS.tex
\section{Krull-Schmidt}\label{sec:Krull-Schmidt}

We now prove a Krull-Schmidt theorem (Theorem \ref{decompuniqness} below)
for decompositions of connected abelian groups into indecomposable subgroups.
This material will not be used in subsequent sections, but provides insight into
the behavior of indecomposable abelian groups.

\begin{lemma}\label{relsupplement}
Let $A < B < C$ be abelian groups of finite Morley rank with $A$ and $B$
definable.  If $A$ has a definable supplement in $C$ then $A$ has
a definable supplement in $B$ and $B$ has a definable supplement in $C$.
\end{lemma}

\begin{proof}
Let $A + A' = C$ with $A' < C$.  Then $A'$ is a supplement for $B$ in $C$.

We claim that $A' \cap B$ is a supplement for $A$ in $B$.
Since $A + A' = C$ and $A' < C$, we have $A \cap A' < A$ and $B \cap A' < B$.
Now $B = (A + A') \cap B = A + (A' \cap B)$.
So $B \cap A'$ is a supplement for $A$ in $B$.
\end{proof}

In particular, the radical $J(A)$ of an abelian group of finite Morley rank
 will not have a supplement in another abelian group containing $A$.

\begin{corollary}\label{relsupplementcor}
Let $B$ be a abelian group of finite Morley rank, and
let $A$ be a subgroup of $B$.  Then $J(A) \leq J(B)$.
\end{corollary} 

\begin{remark}
Let $C$ be $\C^3$ under addition, with predicates for the subgroups
 $A = \C\times1\times1$, $A' = 1\times\C\times1$,
 $B = A + A'$, and $B' = \C \times 1 \times \C$.
Then $A$ has a definable supplement in $B$,
 $B$ has a definable supplement in $C$, and
 $A$ has no definable supplement in $C$.
\end{remark}

Let $A$ be a connected abelian group of finite Morley rank, and 
 let $\mathcal I$ be a family of definable indecomposable subgroups of $G$.
We say $\mathcal I$ is a {\it minimal decomposition} of $A$ into
 indecomposable subgroups if $A = \gen{ \bigcup \mathcal I }$
 and no $X\in \mathcal I$ may be omitted.
We observe that the indecomposable groups in a minimal decomposition
$\mathcal I$ are connected, by Lemma \ref{conind} and the fact that the
group generated by the connected indecomposable groups has finite index in $A$.
By Fact \ref{decomposition}, 
 every connected abelian group of of finite Morley rank has
 a minimal decomposition into indecomposable subgroups.

\begin{lemma}\label{mindecomp_quotient}
Let $A$ be a connected abelian group of finite Morley rank.
Let $\{ A_1,\ldots,A_n\}$ be a minimal decomposition of $A$
 into definable connected indecomposable subgroups.
Then, for any $J \leq A$ with $J^\o \leq J(A)$,  there is  
 a minimal decomposition $\{ A_1 J/J, \ldots, A_n J/J \}$ of $A/J$
 into definable connected indecomposable subgroups.
\end{lemma}

\begin{proof} 
By Fact \ref{Indpushforward},
 $\{ A_1 J/J, \ldots, A_n J/J \}$ is
 a decomposition of $A/J$ into indecomposable subgroups.
Suppose that this decomposition is not minimal,
i.e.\ there is an $i$ such that $A_i \leq J + \sum_{j \neq i} A_j$.
Since $B := \sum_{j \neq i} A_j < A$,
 we find that $B$ and $J$ are supplements to one another.
Since $A$ is connected and $J^\o \leq  J(A)$, 
 $B$ and $J^\o$ are supplements to one another too.
By Lemma \ref{relsupplement},
 $J(A)$ has a supplement in $A$, a contradiction.
\end{proof}

\begin{lemma}\label{directsum}
Let $A$ be a connected abelian group of finite Morley rank.
Let $\{ A_1,\ldots,A_n \}$ be a minimal decomposition of $A$
 into definable connected indecomposable subgroups.
Then there is a subgroup $J \leq A$ such that $J^\o \leq J(A)$
 and $A/J = \bigoplus_{k\leq n} A_k J/J$.
\end{lemma}

\begin{proof}
Let $K := A_n \cap (A_1 + \cdots + A_{n-1})$.
Then $A/K = A_n/K \oplus B$ where
 $B := (A_1 + \cdots + A_{n-1})/K$.
Since $K < A_n$, we have $K^\o \leq J(A_n) \leq J(A)$
 by Corollary \ref{relsupplementcor}.
By Lemma \ref{mindecomp_quotient}, $\{ A_1 K/K, \ldots, A_n K/K \}$
 is a minimal decomposition of $A/K$ into indecomposable subgroups.
So $\{ A_1 K/K, \ldots, A_{n-1} K/K \}$ is
 a minimal decomposition of $B$ into indecomposable subgroups.
By induction, there is a $\tilde{J} \leq B$ such that
 $B/\tilde{J} = \bigoplus_{k < n} (A_k K/K) \tilde{J}/\tilde{J}$,
 and $\tilde{J}^\o \leq J(B)$.

Now let $J$ be the pull-back of $\tilde{J}$
 under the quotient map $x \mapsto x/K$.
Then $A/J = \bigoplus_{k\leq n} A_k J/J$.
If $J^\o$ has a supplement in $A$, then $J^\o$ has a supplement
in $A_1 + \cdots + A_{n-1}$ by Lemma \ref{relsupplement},
 contradicting the fact that $\tilde{J}^\o \leq J(B)$.
\end{proof}

\begin{corollary}\label{directsum_rr}
Let $A$ be a connected abelian group of finite Morley rank, and
let $\{ A_1,\ldots,A_n \}$ be a minimal decomposition of $A$
 into definable connected indecomposable subgroups.
Then $\rr(A) = \sum_{i \leq n} \rr(A_i)$.
\end{corollary}

\begin{proof}
By Lemma \ref{mindecomp_quotient}, $A_1 J(A)/J(A), \ldots, A_n J(A)/J(A)$
 is a minimal decomposition for $A/J(A)$.
By Fact \ref{Indpushforward}, $\rr(A_i J(A)/J(A)) = \rr(A_i)$ for $i\leq n$.
So we may assume that $J(A) = 1$
By Lemma \ref{relsupplement}, $J(A_i) \leq J(A) = 1$ for $i \leq n$.
So the result follows from Lemma \ref{directsum}.
\end{proof}

As one expects, these decompositions are essentially unique, in the
sense of the Krull-Schmidt theorem.

\begin{theorem}[{\cite[\qTheorem 2.40]{BuPhd}}]\label{decompuniqness} 
Let $A$ be a connected abelian group of finite Morley rank.
Also let $\{ A_1,\ldots,A_n \}$ and $\{ B_1,\ldots,B_m \}$
 be two minimal decompositions of $A$ into
 definable connected indecomposable subgroups.
Then the following hold,
 after a suitable reindexing of the set $\{A_1,\ldots,A_n\}$.
\begin{conclusions}
\item $m = n$.
\item $\rr(A_j) = \rr(B_j)$ for all $j \leq n$.
\item For $j \leq n$, we have a minimal decomposition of $A$ given by
$$ A = A_1 + \cdots + A_j + B_{j+1} + \cdots + B_n \mathperiod $$
\end{conclusions}
\end{theorem}

\begin{proof}
We will prove by induction that there is a reindexing such that,
 for every $j \leq n$, we have a minimal decomposition
 $A = A_1 + \cdots + A_{j-1} + B_j + B_{j+1} + \cdots + B_m$
 for $A$, and $\rr(A_j) = \rr(B_j)$.
Suppose inductively that
 $A = A_1 + \cdots + A_{j-1} + B_j + \cdots + B_{m_j}$
 is a minimal decomposition for $A$ with $m_j \leq m$,
 and that $\rr(A_k) = \rr(B_k)$ for $k < j$.
We fix the indices of $A_1,\ldots,A_{j-1}$, and
 may freely reindex the remainder of this list.

By Lemma \ref{directsum},
 there is a subgroup $J \leq A$ such that $J^\o \leq J(A)$
 and $A/J = \bigoplus_{k\leq n} A_k J/J$.
Choose $b\in B_j \setminus K$ where
 $K := J + J(A) + A_1 + \cdots + A_{j-1} + B_{j+1} + \cdots + B_{m_j}$.
There are $a_k \in A_k J$ for $k \leq n$ such that $b = \sum_{k \leq n} a_k$.
By the choice of $b$, there is an $i \geq j$ such that $a_i \notin K$.
Consider the projection $\pi_i : A \to A_i J/J$ given by
 $x \mapsto y_i J/J$ where $x = \sum_{k \leq n} y_k$ with $y_k \in A_k$.
Since $\pi_i(b) = a_i J/J \neq 1$, the map $\pi_i|_{B_j}$ is surjective.
By Fact \ref{Indpushforward}, $\rr(B_j) = \rr(A_i)$,
 concluding part 2 for this stage.

By choice of $b$ and $a_i$, we have
 $(B_j \cap K)^\o \leq J(B_j)$  and $(A_i \cap K)^\o \leq J(A_i)$.
By our inductive assumption and Corollary \ref{directsum_rr},
 $\rr(B_j) = \rk(A/K) = \rk(A) - \rk(K)$.
Now $A_i/J$ is embedded in $A/K$, but $\rr(A_i) = \rr(B_j)$,
so $A = A_1 + \cdots + A_{j-1} + A_i + B_{j+1} + \cdots + B_{m_{j+1}}$,
 with $m_{j+1} = m_j$,
is a decomposition of $A$ into definable indecomposable subgroups.

Since $i \geq j$, we may assume that $i = j$
 by reindexing the $A_k$ with $k \geq j$.
By Corollary \ref{directsum_rr}, $\rr(A)
 = \sum_{k < j} \rr(A_k) + \sum_{k \geq j} \rr(B_k) 
 = \sum_{k \leq j} \rr(A_k) + \sum_{k > j} \rr(B_k)$.
By Corollary \ref{directsum_rr},
 $A = A_1 + \cdots + A_j + B_{j+1} + \cdots + B_m$
 is a minimal decomposition, as desired.

It now follows that $n \leq m$.
By symmetry, $m \leq n$ too.
\end{proof}

% ...
%
%
%
%